\documentclass[12pt,a4paper]{amsart}
\usepackage{graphicx} %(in the preamble) 
\usepackage{preprint}
\begin{document}

\title{Three spheres theorem for $p$-harmonic functions}

\author[V.M. Miklyukov, A. Rasila and M. Vuorinen]
{Vladimir M. Miklyukov, Antti Rasila and Matti Vuorinen}
\address[Vladimir M. Miklyukov]
{Department of Mathematics, Volgograd State University,
Universitetskii prospect 100, Volgograd 400062, RUSSIA,
Fax + tel: +7-8442 471608}\email{miklyuk@mail.ru}\
\address[Antti Rasila]
{Department of Mathematics and Systems Analysis, Aalto University, P.O. Box 1110
0, FI-00076 Aalto, FINLAND, Fax +358-9-451 3016}\
\email{antti.rasila@iki.fi}\
\address[Matti Vuorinen]
{Department of Mathematics, FI-20014 University of Turku,
FINLAND} \email{vuorinen@utu.fi}

%\cgw

%\author{Vladimir M. Miklyukov, Antti Rasila and Matti Vuorinen}
%\email{miklyuk@mail.ru, antti.rasila@tkk.fi, vuorinen@utu.fi}


\thanks{The research was supported by the Academy of Finland, grant 107317 and
by the Foundation Vilho, Yrj\"o ja Kalle V\"ais\"al\"an rahasto of
Finnish Academy of Science and Letters.}

% abstrakti
 \begin{abstract}
   Three spheres theorem type result is proved for the $p$-harmonic 
   functions defined on the complement of $k$-balls in the Euclidean 
   $n$-dimensional space.
 \end{abstract}

\maketitle

\section{Introduction}

A classical theorem by J.\ Hadamard gives the following relation
between the maximum absolute values of an analytic function on
three concentric circles. 

\begin{thm}
Let $R_1<r_1<r_2<r_3<R_2$ and let $f$ be an analytic
function in the annulus $\{z\in \mathbb{C} : R_1<|z|<R_2\}$. Denote the 
maximum of 
$|f(z)|$ on the circle $|z|=r$ by $M(r)$.
Then
$$
M(r_2)^{\log(r_3/r_1)} \leq 
M(r_1)^{\log(r_3/r_2)}M(r_3)^{\log(r_2/r_1)}.
$$
\end{thm}

This result, known as the three circles theorem, was given by
Hadamard without proof in 1896 \cite{Hada}. For a discussion of
the history of this result, see e.g. \cite{Robinson} and 
\cite[pp. 323--325]{Maz}. It is a natural
question, what results of this type can be proved for other classes
of functions. For example, a version of Hadamard's theorem can be proved for
subharmonic functions in $\Rn$, $n\geq 2$, see \cite[pp. 128--131]{Protter}.

%\section{Preliminaries}

Some generalizations of the three circles theorem will be studied here. 
For the formulation of our main result, Theorem \ref{main}, we recall
some standard notation and definitions from the book \cite{HKM}. We 
will consider solutions $v\colon \Omega \to \R$ of the 
{\it $p$-Laplace equation}
\begin{equation}
\label{star1}
\diver\big(|\nabla v|^{p-2}\nabla v\big)=0, \qquad 1<p<\infty,
\end{equation}
on an open set $\Omega\subset \Rn$ in the sense that will be described 
shortly. When $p=2$ equation (\ref{star1})
reduces to the Laplace equation $\Delta u=0$, whose solutions,
harmonic functions, are studied in the classical potential theory.  When 
$p\neq 2$ equation (\ref{star1}) is nonlinear and degenerates 
at the zeros of the gradient of $v$. It follows that the solutions,
{\it $p$-harmonic functions}, need not be in $C^2(\Omega)$ and the
equation must be understood in the weak sense. A 
weak solution of (\ref{star1}) is a function $v$ in the Sobolev space 
$W^{1,p}_{\loc}(\Omega)$ such that
\begin{equation}
\label{star1weak}
\int_{\Omega}\langle |\nabla v|^{p-2} \nabla v,\nabla 
\varphi\rangle\,dm=0
\end{equation}
for all $\varphi\in C_0^\infty(\Omega)$, where 
$\langle\cdot,\cdot\rangle$ 
denotes the scalar product of vectors in $\Rn$, and $m$ is the
Lebesgue measure in $\Rn$.

It is easy to see that for all $\varphi\in 
C_0^\infty(\Omega)$ and $v\in C^2(\Omega)$,
$$
\int_{\Omega}\langle |\nabla v|^{p-2} \nabla v,\nabla 
\varphi\rangle\,dm=-\int_{\Omega}\varphi\,{\rm div}\,(|\nabla v|^{p-2} 
\nabla v)\,dm
$$
and, consequently, each $C^2$-solution to (\ref{star1}) is a weak solution to 
(\ref{star1}).

Fix an integer $k$, $1\le k\le n$ and a real number $t\geq 0$. The sets
$B_k(t)=\{x\in \Rn:d_k(x)<t\}$ 
and $\Sigma_k(t)=\{x\in \Rn:d_k(x)=t\}=\partial B_k(t)$,
 where 
$d_k(x)=\Bigl(\sum\limits_{i=1}^kx_i^2\Bigr)^{1/2}$,
are respectively called $k$-ball and $k$-sphere in $\Rn$. For $k=n$ the 
$k$-ball $B_k(t)$ coincides with the standard Euclidean ball $B^n(t)$ 
and the $k$-sphere $\Sigma_k(t)$ is the Euclidean sphere $S^{n-1}(t)$.
In particular, the symbol $\Sigma_k(0)$ below denotes the $k$-sphere 
with the radius $0$, i.e.
$$
\Sigma_k(0)=\{x=(x_1,\ldots,x_k,\ldots,x_n):x_1=\ldots=x_k=0\}.
$$

Let $0<\alpha<\beta<\infty$ be fixed and let 
$$
D_{\alpha,\beta}=\{x\in \Rn: \alpha<d_k(x)<\beta\}.
$$ 
For $k=1$ the set $D_{\alpha,\beta}$ is the union of the two layers
between two parallel hyperplanes. For $1<k<n$ the boundary of
the domain $D_{\alpha,\beta}$ consists of two coaxial cylindrical 
surfaces.

\begin{figure}[t]
\begin{center}
\includegraphics[width=5.5cm]{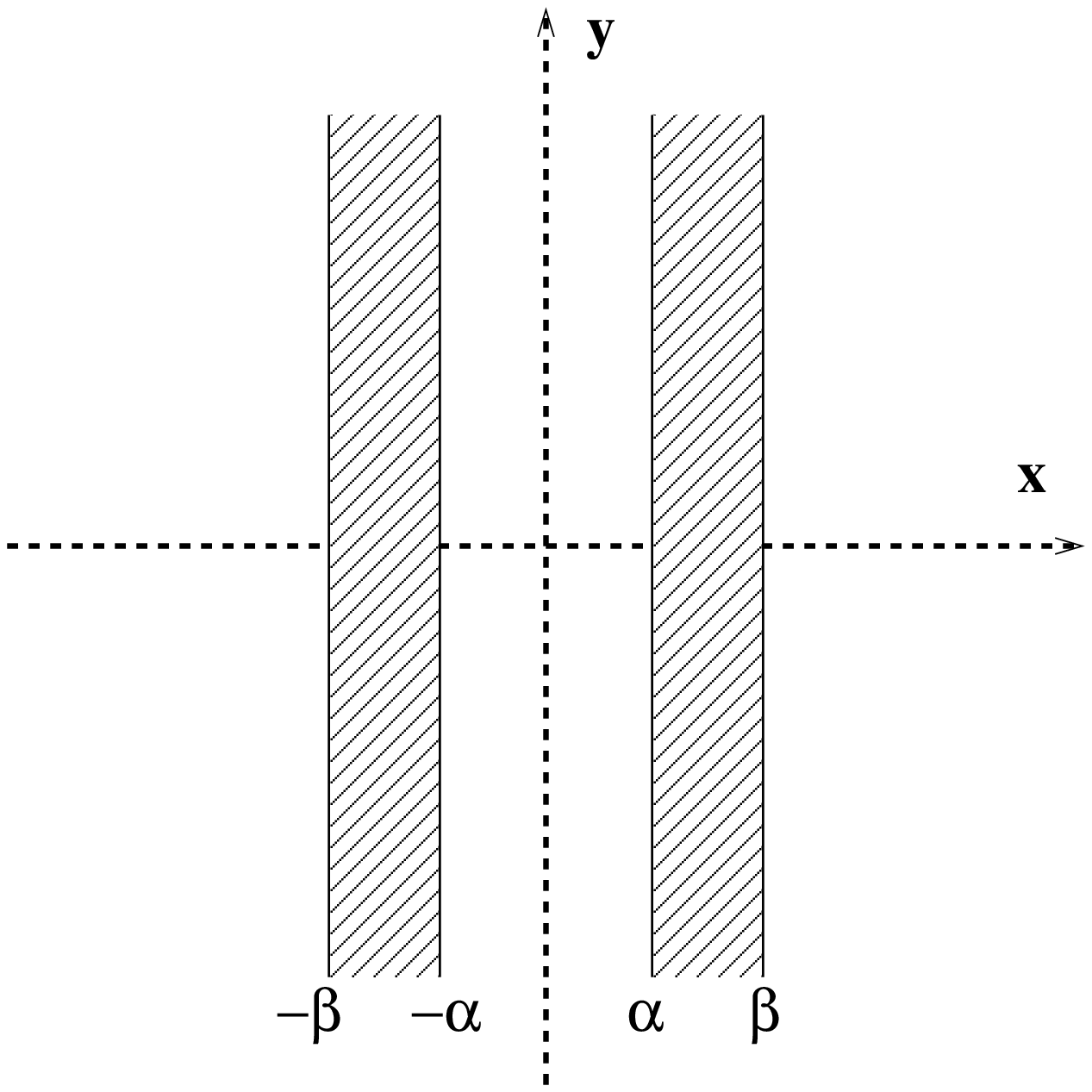}
\qquad
\includegraphics[width=5.5cm]{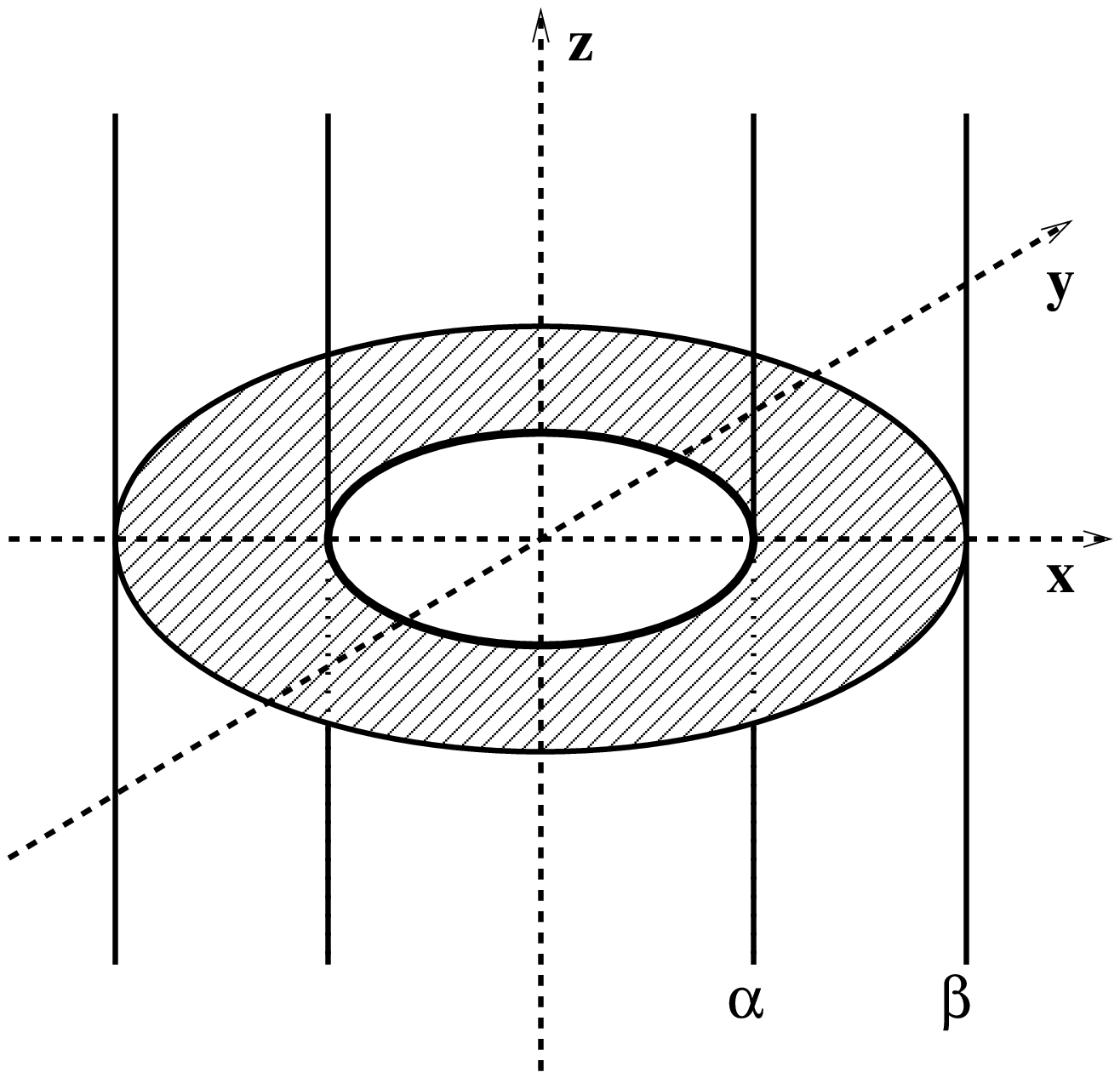}

\caption{$1$-annulus $D_{\alpha,\beta}$ in $\R^2$ (left)
and $2$-annulus $D_{\alpha,\beta}$ in $\R^3$ (right).}
\label{fig:fig1}
\end{center}
\end{figure}

Let $v\in C^0(D_{r,R})$, and let $M(r)=\limsup_{z\to\Sigma_k(r)}v(z)$. 
Suppose that $M(R) > M(r)$. Consider the function
$$
v_{r,R}(x)=
\frac{v(x)-M(r)}{M(R)-M(r)},
$$
for $r<R$. Clearly, $\limsup_{z\to\Sigma_k(r)}v_{r,R}(z)
\leq 0$ and 
$\limsup_{z\to\Sigma_k(R)}v_{r,R}(z) \leq 1$.
Let 
$$
\xi(r,t)=\int_r^t s^{(1-k)/(p-1)}ds,
\text{ and }
u_0^{k,p}(t)=\frac{\xi(r,t)}{\xi(r,R)}.
$$
Let $u(x)=u_0^{k,p}\big(d_k(x)\big)$ for $x\in D_{r,R}$.
It is clear (see Lemma \ref{psolu}) that $u$ is a $C^2$-solution to 
(\ref{star1}).
We have                                           
$$
u(x)|_{\Sigma_k(r)}\equiv 0,\,\,\,
u(x)|_{\Sigma_k(R)}\equiv 1,
$$
and
\begin{equation}
\label{star2}
u(x)\geq v_{r,R}(x)\textrm{ if }
x\in\Sigma_k(r)\textrm{ or }
x\in\Sigma_k(R).
\end{equation}

\section{Main results}

We will prove the following Hadamard type theorem for the $p$-harmonic
functions defined on the complement of a $k$-ball. We use the method of
proof from \cite{MV}.

\begin{thm}
\label{main}
Let $1<p<\infty$, $R>r>0$ and let
$v(x)\in W^{1,p}_{\loc}(D_{r,\infty})$
be a 
continuous weak solution of (\ref{star1}) such that
\begin{equation}
\label{star4}
\int_{r}^{\infty}dt
\bigg(\int_{\Sigma_k(t)}|v_{r,R} -u|^2
\big(|\nabla v_{r,R}|^{|p-2|}+|\nabla u|^{|p-2|}\big)
d{\Hm}^{n-1}\bigg)^{-1} = \infty,
\end{equation}
where $\Hm^{n-1}$ is the $(n-1)$-dimensional Hausdorff measure.
Then for all $t\in(r,R)$,
\begin{equation}
\label{star3}
M(t) \leq \big(M(R)-M(r)\big)u_0^{k,p}(t)+M(r).
\end{equation}
\end{thm}

Note that for $k=n$ (\ref{star3}) follows immediately from the 
comparison principle, see \cite[p. 133]{HKM}.

\begin{cor}
\label{cora}
Let $1<p<\infty$, $R>r>0$ and let
$v(x)\in W^{1,p}_{\loc}(D_{r,\infty})$, be a 
continuous weak solution of (\ref{star1}) such that
\begin{equation}
\label{star4b}
\lim_{S\to\infty}
\frac{1}{S^2}
\int_{D_{r,S}}|v_{r,R} -u|^2
\big(|\nabla v_{r,R}|^{|p-2|}+|\nabla u|^{|p-2|}\big)
dm = 0.
\end{equation}
Then for all $t\in(r,R)$ the inequality (\ref{star3}) holds.
\end{cor}

\begin{cor}
\label{corb}
Let $1<p<\infty$, $R>r>0$ and let
$v(x)\in W^{1,p}_{\loc}(D_{r,\infty})$
be a 
continuous weak solution of (\ref{star1}) such that
$$
\int_{D_{r,\infty}}|v_{r,R} -u|^2
\big(|\nabla v_{r,R}|^{|p-2|}+|\nabla u|^{|p-2|}\big)
dm \leq M < \infty.
$$
Then for all $t\in(r,R)$ the inequality (\ref{star3}) holds.
\end{cor}

%\subsection*{Variational kernels}

For the formulation of a result of S.~Granlund \cite{Gran}, 
Theorem \ref{granthm} below, we introduce some notation and
terminology. Let $p>1$, $\Omega \subset \Rn$ be a bounded domain and let
$F\colon \Omega \times \Rn \to \R$ be such that the following conditions hold.
\begin{enumerate}
\item
There are constants $\beta>\alpha>0$ such that for a.e.
$x\in \Omega$
$$
\alpha|z|^p\leq F(x,z) \leq \beta|z|^p.
$$
\item
For a.e. $x\in \Omega$ the function $z\mapsto F(x,z)$ is convex.
\item
The function $x\mapsto F(x,\nabla u(x))$ is measurable for all
$u\in W^{1,p}(\Omega)$.
\end{enumerate}
Let
$$
I(u)= \int_\Omega F\big(x,\nabla u(x)\big)\, dm.
$$
A function $u\in W^{1,p}(\Omega)$ is a {\it subminimum} in $\Omega$ if 
$I(u)\leq I(u-\eta)$ for all non-negative $\eta \in W^{1,p}_0(\Omega)$.
Let 
$$
M(r)=\esssup_{x\in \overline{B}^n(r)}u(x),\qquad
\overline{B}^n(r)\subset \Omega.
$$
The following Hadamard type theorem was proved by S.~Granlund
in \cite{Gran}. 

\begin{thm}
\label{granthm}
Let $u$ be a subminimum of
$$
I(u)= \int_\Omega F\big(x,\nabla u(x)\big) dm,
$$
$r_1<r<r_2$, and $\overline{B}^n(r_2)\subset \Omega$. Then $u$
is bounded from above, and there is a constant
$$\lambda=\lambda(n,p,r/r_1,r_2/r,\alpha/\beta),$$ $0<\lambda<1$
such that
$$
M(r)\leq \lambda M(r_1) + (1-\lambda) M(r_2).
$$
\end{thm}

Since $p$-harmonic functions minimize (see e.g. \cite[p. 59]{HKM}) the 
integral 
$$
I(u)=\int_\Omega|\nabla u|^p\,dm,
$$
Theorem \ref{granthm} is related to Theorem \ref{main} with 
$k=n$.

\section{Preliminaries}

We start by recalling some basic properties of the Sobolev spaces
from \cite{HKM}. Let $\Omega$ be a nonempty open set in $\Rn$.

\begin{lem}
\label{sobprod}
\cite[Theorem 1.24]{HKM}
Let $u\in W^{1,p}_0(\Omega)$ and $v \in W^{1,p}(\Omega)$ be bounded. Then
$uv\in W^{1,p}_0(\Omega)$.
\end{lem}

\begin{lem}
\label{appro}
\cite[Lemma 3.11]{HKM}
If $v\in W^{1,p}(\Omega)$ is a weak solution of (\ref{star1}) in 
$\Omega$, then
$$
\int_{\Omega}\langle |\nabla v|^{p-2} \nabla v,\nabla 
\varphi\rangle\,dm=0
$$
for all $\varphi\in W^{1,p}_0(\Omega)$.
\end{lem}

%We will also make use of Federer's co-area formula formulated below.

\begin{thm}
\cite[p. 99]{Evans}
\label{coarea}
Let $f\colon{\R}^n\to \R$ be a locally Lipschitz mapping.
Let $E\subset{\R}^n$ be an $n$-measurable set and $g\colon E\to{\R}$
be a nonnegative measurable function. Then
\begin{equation}
\label{2eqare}
\int\limits_E g(x)|\nabla f(x)|\,dx_1\cdots dx_n
=\int\limits_{\R}
\bigg(\sum_{x\in f^{-1}(y)}g(x)\bigg)
\,d{\Hm}^n (y).
\end{equation}
%where $\left|J(x,f)\right|$ is defined as in \cite[p. 91]{Evans}.
\end{thm}

%\subsection*{Two elementary lemmas}

\begin{lem}
\label{psolu}
Let $1<p<\infty$, $0<r<d_k(x)$ and fix an integer $1\leq k \leq n$. Then
$$
u(x) = \int_r^{d_k(x)}s^{\frac{1-k}{p-1}}\,ds
$$
is a solution of (\ref{star1}), i.e.
$$
\sum_{i=1}^n\Big\{\frac{\partial}{\partial x_i}\big(u_{x_i}[u_{x_1}^2+
\ldots+u_{x_n}^2]^{\frac{p-2}{2}}\big)\Big\}=0.
$$
\end{lem}

\proof
We note that
$$
\frac{\partial}{\partial x_i} d_k(x)=\frac{x_i}{d_k(x)},
$$
and hence $u_{x_i}=x_i d_k(x)^{\frac{1-k}{p-1}-1}$. Then
\begin{eqnarray*}
u_{x_i}\big(u_{x_1}^2+...+u^2_{x_k}\big)^{\frac{p-2}{2}}
&=&
x_id_k(x)^{\frac{1-k}{p-1}-1}\Big[d_k(x)^{\frac{2(1-k)}{p-1}-2}
\big(\sum_{j=1}^k x_j^2\big)\Big]^{\frac{p-2}{2}}\\
&=&
x_id_k(x)^{\frac{1-k}{p-1}-1}d_k(x)^{\frac{(1-k)(p-2)}{p-1}}
= x_id_k(x)^{-k}.
\end{eqnarray*}
It follows that
\begin{eqnarray*}
\sum_{i=1}^k\frac{\partial}{\partial x_i}\big(x_id_k(x)^{-k}\big)
&=&
\sum_{i=1}^k d_k(x)^{-k}
-k \sum_{i=1}^k x_i^2 d_k(x)^{-k-2}\\
&=&
kd_k(x)^{-k} - kd_k(x)^{-k-2}\big(\sum_{i=1}^k x_i^2\big)
=0.
\end{eqnarray*}
\qed

Next we will prove two lemmas which are used later in the proof
of Theorem \ref{main}.

\begin{lem}
\label{elem1}
Let $a > b>0$, $p>1$. Then
%It is easy to see that
\begin{equation}
\label{eq19}
C_1\,{\frac{a^{p-1}-b^{p-1}}{a-b}}\le{\frac{a^{p-1}+b^{p-1}}{a+b}}
\le C_2\,{\frac{a^{p-1}-b^{p-1}}{a-b}},
%\quad a\ge b>0,\quad p>1,
\end{equation}
with some constants $C_1,\, C_2>0$.
\end{lem}

\proof
We examine the function 
\begin{equation}
\label{limi}
g_1(x)={\frac{(x^{p-1}+1)(x-1)}{(x^{p-1}-1)(x+1)}},\qquad
x> 1. $$ It is clear that $$ \lim_{x\to 1}g_1(x)={\frac{1}{p-1}},\qquad
\lim_{x\to
\infty}g_1(x)=1. 
\end{equation}
It is sufficient to find positive bounds for $g_1(x)$ for $x>1$. 
We will prove that the bounds are in fact given by (\ref{limi}).
First we note that
$$
\left\{\begin{array}{ll}
(p-2)(x^p-1)+p(x-x^{p-1})
<  0, &\text{ for }p\in(1,2),\\
(p-2)(x^p-1)+p(x-x^{p-1})
=  0, &\text{ for }p=2,\\
(p-2)(x^p-1)+p(x-x^{p-1})
>  0, &\text{ for }p>2,
\end{array}\right.
$$
and
$$
\left\{\begin{array}{ll}
x-x^{p-1}
<  0, &\text{ for }p\in(1,2),\\
x-x^{p-1}
=  0, &\text{ for }p=2,\\
x-x^{p-1}
>  0, &\text{ for }p>2.
\end{array}\right.
$$
Hence
$$
\left\{\begin{array}{ll}
g_1(x)
\in  \big(1,1/(p-1)\big), &\text{ for }p\in(1,2),\\
g_1(x) 
=  1, &\text{ for }p=2,\\
g_1(x)
\in \big(1/(p-1),1\big), &\text{ for }p>2.
\end{array}\right.
$$
\qed

\begin{lem}
Let $a > b>0$. Then
\begin{equation}
\label{eq20}
C_3\,(a^{p-2}+b^{p-2}) \le{\frac{a^{p-1}-b^{p-1}}{a-b}} \le
C_4\,(a^{p-2}+b^{p-2}),
\end{equation}
for $p\ge 2$, and 
\begin{equation}
\label{eq20'}
C_3\,(a^{2-p}+b^{2-p})^{-1} \le{\frac{a^{p-1}-b^{p-1}}{a-b}} \le
C_4\,(a^{2-p}+b^{2-p})^{-1},
\end{equation}
for $p\in(1,2]$ with some constants $C_3,\, C_4>0$.
\end{lem}

\proof
The proof is similar to that of Lemma \ref{elem1}. First we
study the function
$$
g_2(x)=\frac{x^{p-1}-1}{(x-1)(x^{p-2}+1)}.
$$
As in Lemma \ref{elem1}, it is sufficient for (\ref{eq20}) to find
positive bounds for $g_2(x)$ for $x>0$. We note that
$\lim_{x\to 1}g_2(x)=(p-1)/2$ and $\lim_{x\to\infty}g_2(x)=1$. 
We obtain
$$
\left\{\begin{array}{ll}
(p-3)(1-x^{p-1})+(p-1)x(1-x^{p-3})
<  0, &\text{ for }p\in(1,3),\\
(p-3)(1-x^{p-1})+(p-1)x(1-x^{p-3})
=  0, &\text{ for }p=3,\\
(p-3)(1-x^{p-1})+(p-1)x(1-x^{p-3})
>  0, &\text{ for }p>3,
\end{array}\right.
$$
and
$$
\left\{\begin{array}{ll}
x(x^{p-3}-1)
<  0, &\text{ for }p\in(1,3),\\
x(x^{p-3}-1)
=  0, &\text{ for }p=3,\\
x(x^{p-3}-1)
>  0, &\text{ for }p>3.
\end{array}\right.
$$
It follows that
$$
\left\{\begin{array}{ll}
g_2(x)
\in  \big((p-1)/2,1\big), &\text{ for }p\in(1,3),\\
g_2(x) 
=  1, &\text{ for }p=3,\\
g_2(x)
\in \big(1,(p-1)/2\big), &\text{ for }p>3.
\end{array}\right.
$$
To prove (\ref{eq20'}) we study the function
$$
g_3(x)=\frac{(x^{p-1}-1)(x^{2-p}+1)}{x-1}.
$$
Now $\lim_{x\to 1}g_3(x)=2(p-1)$ and $\lim_{x\to\infty}g_3(x)=1$.
Again, we have
$$
\left\{\begin{array}{ll}
(-2p+3)(x-1)+(x^{p-1}-x^{2-p})
<  0, &\text{ for }p\in(1,3/2),\\
(-2p+3)(x-1)+(x^{p-1}-x^{2-p})
=  0, &\text{ for }p=3/2,\\
(-2p+3)(x-1)+(x^{p-1}-x^{2-p})
>  0, &\text{ for }p>3/2,
\end{array}\right.
$$
and
$$
\left\{\begin{array}{ll}
x^{p-1}-x^{2-p}
<  0, &\text{ for }p\in(1,3/2),\\
x^{p-1}-x^{2-p}
=  0, &\text{ for }p=3/2,\\
x^{p-1}-x^{2-p}
>  0, &\text{ for }p>3/2,
\end{array}\right.
$$
and thus
$$
\left\{\begin{array}{ll}
g_3(x)
\in  \big(2(p-1),1\big), &\text{ for }p\in(1,3/2),\\
g_3(x) 
=  1, &\text{ for }p=3/2,\\
g_3(x)
\in \big(1,2(p-1)\big), &\text{ for }p>3/2.
\end{array}\right.
$$
\qed

\section{Proof of Theorem \ref{main}}

Suppose the contrary, that is, there exists $x_0\in D_{r,R}$ such
that
\begin{equation}
v(x_0) > \big(M(R)-M(r)\big)u(x_0)+M(r),
\end{equation}
or
$$
v_{r,R}(x_0) > u(x_0).
$$
Fix some $\varepsilon_0>0$, for which
$$
v_{r,R}(x_0) -u(x_0) > \varepsilon_0.
$$
Consider the set
$$
U=\{x\in D_{r,R} : v_{r,R}(x)-u(x)>\varepsilon_0\}\neq\emptyset.
$$
%This set is not empty. 
Choose a component $O$ of $U$ such that $x_0\in O$.
It is clear that $\overline{O}\cap\partial D_{r,R}=\emptyset$ and
$\big(v_{r,R}(x)-u(x)\big)|_{\partial O} =0$.
Fix $\varepsilon_2 > \varepsilon_1>0$ and the balls
$O_1=B_k(x_0,\varepsilon_1)$, $O_2=B_k(x_0,\varepsilon_2)$. Let
$\varphi(x)=\eta(d_k(x))$ be a locally
Lipschitz function with the properties:
\begin{equation}
\label{nb1}
\left\{\begin{array}{lll}
\varphi\equiv 1 &\textrm{for all }x\in O_1,\\
\varphi\equiv 0 &
\textrm{for all }x\in D_{r,R}\setminus O_2.
\end{array}\right.
\end{equation}
Then the function $\psi = \big(v_{r,R}(x)-u(x)\big)\varphi^2$
has a support $\supp \psi \subset \overline{O}_2$ and by Lemma 
\ref{sobprod} $\psi \in W^{1,p}_0(\Omega)$ for all $\Omega$ such that
$\supp \psi \subset \Omega$. Since $v_{r,R}$ and
$u$ are generalized solutions of (\ref{star1})
we have by Lemma \ref{appro} 
\begin{multline*}
\int_{\supp \psi}
\langle \nabla \psi,|\nabla v_{r,R}|^{p-2}\nabla v_{r,R}
- |\nabla u|^{p-2}\nabla u\rangle\, dm\\
=
\int_{\supp \psi}
\langle \nabla \psi,|\nabla v_{r,R}|^{p-2}\nabla v_{r,R} \rangle\,dm
-\int_{\supp \psi} \langle \nabla \psi, |\nabla u|^{p-2}\nabla u \rangle\,dm
=0.
\end{multline*}
%where $dm=dx_1\ldots dx_n$.
Next, we note that
$$
\nabla \psi=
\varphi^2(\nabla v_{r,R} -\nabla u)
+2\varphi (v_{r,R}-u)\nabla\varphi.
$$
Thus, we may write
\begin{eqnarray*}
0 &=&
\int_{\supp \psi}\langle \nabla\psi,|\nabla v_{r,R}|^{p-2}
\nabla v_{r,R}-|\nabla u|^{p-2}\nabla u\rangle\, dm \\
&=&
\int_{O\cap O_2}
\langle \nabla \psi,|\nabla v_{r,R}|^{p-2}\nabla v_{r,R}
- |\nabla u|^{p-2}\nabla u\rangle\, dm\\
&=&
\int_{O\cap O_2}
\varphi^2\langle \nabla v_{r,R}-\nabla u,
|\nabla v_{r,R}|^{p-2}\nabla v_{r,R}
- |\nabla u|^{p-2}\nabla u\rangle\, dm\\
&&+
2\int_{O\cap O_2}
\varphi(v_{r,R}-u)\langle \nabla \varphi,
|\nabla v_{r,R}|^{p-2}\nabla v_{r,R}
- |\nabla u|^{p-2}\nabla u\rangle\, dm
\end{eqnarray*}
or
\begin{multline*}
\int_{O\cap O_2}\varphi^2\langle
\nabla v_{r,R}-\nabla u,|\nabla v_{r,R}|^{p-2}\nabla v_{r,R}
-|\nabla u|^{p-2}\nabla u\rangle\, dm\\
=
-2 \int_{O\cap O_2}
\varphi(v_{r,R}-u)\langle\nabla \varphi,
|\nabla v_{r,R}|^{p-2}\nabla v_{r,R}-|\nabla u|^{p-2}
\nabla u\rangle\, dm
\end{multline*}
or
\begin{multline}
\label{nb2}
\bigg|\int_{O\cap O_2}
\varphi^2\langle\nabla v_{r,R}-\nabla u,
|\nabla v_{r,R}|^{p-2}\nabla v_{r,R}-
|\nabla u|^{p-2}\nabla u\rangle\, dm\bigg|\\
\leq 
2\int_{O\cap O_2}
|\varphi||v_{r,R}-u||\nabla \varphi|
\big||\nabla v_{r,R}|^{p-2}\nabla v_{r,R}-
|\nabla u|^{p-2}\nabla u\big| dm.
\end{multline}
%In the same way as inequality (2.4), (2.5) in \cite{MV} we
%prove that at every point, where $v_{r,R}$ has a 
%differential, the following inequalities hold:
Let 
$$
\Phi(\lambda)=|\nabla(\lambda v_{r,R}+(1-\lambda)u)|^{p-2}
\nabla(\lambda
v_{r,R}+(1-\lambda)u)
$$ for $\lambda\in [0,1]$, and note that $$ \Phi(0)=|\nabla
u|^{p-2}\nabla
u\quad\mbox{and}\quad \Phi(1)=|\nabla v_{r,R}|^{p-2}\nabla v_{r,R}.
$$
Now we write
\begin{multline}
\label{eq8} 
|\nabla
v_{r,R}|^{p-2}\nabla v_{r,R} - |\nabla u|^{p-2}\nabla u = \Phi(1)
-
\Phi(0)=\int\limits_0^1 \Phi'(\lambda)\,d\lambda\\
=\int\limits_0^1 \big[(\nabla v_{r,R} - \nabla u)\, |\nabla
(\lambda v_{r,R} +(1-\lambda)u)|^{p-2}
 +(p-2)\nabla (\lambda v_{r,R}+(1-\lambda)u)\,\\
\cdot|\nabla (\lambda v_{r,R}+(1-\lambda)u)|^{p-4}\,\langle\nabla
v_{r,R}-\nabla
u, \nabla (\lambda v_{r,R}+(1-\lambda)u)\rangle\big]\,d\lambda,
\end{multline}
and obtain
\begin{multline}
\label{eq9} 
 \langle\nabla v_{r,R}-\nabla
u,\,|\nabla v_{r,R}|^{p-2}\nabla v_{r,R} - |\nabla u|^{p-2}\nabla
u\rangle\\
%\begin{equation}
=|\nabla v_{r,R}-\nabla u|^2\int\limits_0^1 |\nabla (\lambda
v_{r,R} + (1-\lambda)u)|^{p-2}\,d\lambda \\
+(p-2)\int\limits_0^1 |\nabla (\lambda v_{r,R} +
(1-\lambda)u)|^{p-4} \langle\nabla v_{r,R}-\nabla u,\,\nabla
(\lambda
v_{r,R} + (1-\lambda)u) \rangle^2\,d\lambda.
\end{multline}
If $p\geq 2$ then
\begin{multline}
\label{nb3a}
\langle\nabla v_{r,R}-\nabla u,
|\nabla v_{r,R}|^{p-2}\nabla v_{r,R}-
|\nabla u|^{p-2}\nabla u\rangle\\
\geq
|\nabla v_{r,R}-\nabla u|^2
\int_0^1\big|\nabla\big(\lambda v_{r,R}
+(1-\lambda)u\big)\big|^{p-2}d\lambda.
\end{multline}
If $p<2$, we have
\begin{multline*}
|\nabla v_{r,R}-\nabla u|^2\int_0^1
\big|\nabla\big(\lambda v_{r,R}
+(1-\lambda)u\big)\big|^{p-2}d\lambda\\
+(p-2)\int_0^1
\big|\nabla\big(\lambda v_{r,R}+(1-\lambda)u\big)\big|^{p-4}
\langle\nabla v_{r,R}-\nabla u,
\nabla\big(\lambda v_{r,R}+(1-\lambda)u\big)\rangle^2\, d\lambda\\
\geq (p-1)|\nabla v_{r,R}-\nabla u|^2\int_0^1
\big|\nabla\big(\lambda v_{r,R}+(1-\lambda)u\big)\big|^{p-2}d\lambda.
\end{multline*}
%Again, in the same way as \cite[2.5]{MV}, we obtain
This together with (\ref{eq9}) gives
\begin{multline}
\label{nb3b}
\langle\nabla v_{r,R}-\nabla u,
|\nabla v_{r,R}|^{p-2}\nabla v_{r,R}-|\nabla u|^{p-2}\nabla u\rangle
\\
\geq (p-1)|\nabla v_{r,R}-\nabla u|^2\int_0^1
\big|\nabla\big(\lambda v_{r,R}+(1-\lambda)u\big)\big|^{p-2}d\lambda,
\,\,\,1<p\leq 2.
\end{multline}
%As in \cite[2.6]{MV}, 
It follows from (\ref{eq8}) that for every $p>1$, 
\begin{multline}
\label{nb4}
\big||\nabla v_{r,R}|^{p-2}\nabla v_{r,R}-
|\nabla u|^{p-2}\nabla u\big|
\leq
C_5|\nabla v_{r,R}-\nabla u|
\int_0^1\big|\nabla\big(\lambda v_{r,R}
+(1-\lambda)u\big)\big|^{p-2}d\lambda,
\end{multline}
at every point where
$v_{r,R}$ has differential. Here $C_5=1+|p-2|$. 
Setting 
$$
I(p) = \int_0^1\big|\nabla(\lambda v_{r,R}
+(1-\lambda)u\big|^{p-2} d\lambda
$$
and using (\ref{nb2}), (\ref{nb3a}), (\ref{nb3b})
and (\ref{nb4}) we obtain
\begin{equation}
\label{nb5}
\int_{O\cap O_2}
\varphi^2 I(p)|\nabla v_{r,R} -\nabla u|^2 dm
\leq C_6\int_{O\cap O_2}
I(p)|\varphi||v_{r,R} -u||\nabla\varphi|
|\nabla v_{r,R}- \nabla u|dm,
\end{equation}
where $C_6=2C_5/\min\{1,p-1\}$.

%%%%%%%%%%%%%%%%

We note that 
\begin{equation*}
 |\nabla (\lambda v_{r,R} +
(1-\lambda)u)|^2=\\
\lambda^2\,|\nabla
v_{r,R}|^2+2\lambda(1-\lambda)\langle\nabla v_{r,R},\,\nabla u\rangle + (1-
\lambda)^2\,|\nabla
u|^2, 
\end{equation*}
and therefore
\begin{equation}
\label{eq13}
\big|\lambda |\nabla v_{r,R}| - (1-\lambda)\,|\nabla u|\big|\\ \le
|\nabla (\lambda v_{r,R} + (1-\lambda)u)|\le \lambda |\nabla v_{r,R}| +
(1-\lambda)|\nabla u|
\end{equation}
for an arbitrary $\lambda\in  [0,1]$.
Let $p\ge 2$. We suppose that $|\nabla v_{r,R}| > |\nabla u|$. Then by
(\ref{eq13}),
\begin{multline}
I(p)\le\int\limits_0^1 \left(\lambda (|\nabla v_{r,R}|-|\nabla u|) +
|\nabla
u|\right)^{p-2}\,d\lambda\\={\frac{1}{|\nabla v_{r,R}|-|\nabla
u|}}\int\limits_{|\nabla
u|}^{|\nabla v_{r,R}|} s^{p-2}\,ds
={\frac{1}{p-1}}{\frac{|\nabla v_{r,R}|^{p-1}-|\nabla
u|^{p-1}}{|\nabla v_{r,R}|-|\nabla u|}}.
\label{eq14}
%\end{equation}
\end{multline}
Next by (\ref{eq13}),
\begin{eqnarray*}
I(p)&\ge&\int\limits_0^1 \big|\lambda |\nabla v_{r,R}|
-(1-\lambda)\,|\nabla u| \big|^{p-2}\,d\lambda\\
&=&\int\limits_0^1\big|\lambda (|\nabla
v_{r,R}|+|\nabla u|)-|\nabla u| \big|^{p-2}\,d\lambda \\
 &=&\int\limits_s^1\left(\lambda
(|\nabla v_{r,R}|+|\nabla u|)-|\nabla u| \right)^{p-2}\,d\lambda\\
&&+\int\limits_0^s
\left(|\nabla u|-\lambda (|\nabla v_{r,R}|+|\nabla u|)
\right)^{p-2}\,d\lambda,\nonumber
\end{eqnarray*}
 where
\begin{equation}
\label{eq15}
s={\frac{|\nabla u|}{|\nabla v_{r,R}|+|\nabla u|}}.
\end{equation}
By computing both of the last two integrals, we obtain
\begin{equation}
\label{eq16}
I(p)\ge{\frac{1}{p-1}}{\frac{|\nabla v_{r,R}|^{p-1}+|\nabla
u|^{p-1}}{|\nabla v_{r,R}|+|\nabla u|}}.
\end{equation}

Let $1<p<2$. As above, we assume $|\nabla v_{r,R}| > |\nabla u|$. Then by
(\ref{eq13}),
\begin{eqnarray*}
I(p)&\le&\int\limits_0^1{{}{\big|\lambda |\nabla v_{r,R}| 
-(1-\lambda)\,
|\nabla
u|\big|^{2-p}}}\,d\lambda\\
&=&\int\limits_0^1{{}{\big|\lambda(|\nabla
v_{r,R}|+|\nabla
u|)- |\nabla u|\big|^{2-p}}}\,d\lambda\\
 &=&\int\limits_0^s\big(|\nabla u|-\lambda
(|\nabla v_{r,R}|+|\nabla u|) \big)^{2-p}\,d\lambda\\
&&
+\int\limits_s^1\big(\lambda(|\nabla
v_{r,R}|+|\nabla u|)-|\nabla u| \big)^{2-p}\,d\lambda,%\nonumber
\end{eqnarray*}
where $s$ is defined in
(\ref{eq15}), and hence	
\begin{equation}
\label{eq17}
I(p)\le {\frac{1}{p-1}}{\frac{|\nabla v_{r,R}|^{p-1}+|\nabla
u|^{p-1}}{|\nabla v_{r,R}|+|\nabla u|}}.
\end{equation}

By (\ref{eq14}), it follows that
\begin{equation}
\label{eq18}
I(p)\ge{\frac{1}{p-1}}{\frac{|\nabla v_{r,R}|^{p-1}-|\nabla
u|^{p-1}}{|\nabla v_{r,R}|-|\nabla u|}}.
\end{equation}

Setting $a=|\nabla v_{r,R}|$ and $b=|\nabla u|$ in (\ref{eq19}),
(\ref{eq20}) and
(\ref{eq20'}), we can obtain by (\ref{eq14}), (\ref{eq16}), (\ref{eq17})
and (\ref{eq18}), for $p\ge 2$
\begin{equation}
\label{nb6a}
C_7\,\left(|\nabla v_{r,R}|^{p-2}+|\nabla u|^{p-2}\right)\le I(p)\le
C_8\,\left(|\nabla v_{r,R}|^{p-2}+|\nabla u|^{p-2}\right),
\end{equation}
or
\begin{equation}
\label{nb6b}
C_7\,\left(|\nabla v_{r,R}|^{2-p}+|\nabla u|^{2-p}\right)^{-1}\le I(p)
\le
C_8\,\left(|\nabla v_{r,R}|^{2-p}+|\nabla u|^{2-p}\right)^{-1},
\end{equation}
$1<p\le 2$,
with some constants $C_7,\,C_8>0$.  The case ${|\nabla v_{r,R}|<|\nabla
u|}$ is analogous.
%
%%%%%%%%
%
%From the inequalities (2.16), (2.17) in \cite{MV} it follows that
%\begin{equation}
%\label{nb6a}
%C_3\big(|\nabla v_{r,R}|^{p-2} - |\nabla u|^{p-2}\big)
%\leq I(p)\leq C_4\big(|\nabla v_{r,R}|^{p-2}+|\nabla u|^{p-2}\big),
%\end{equation}
%for $p\geq 2$ and
%\begin{equation}
%\label{nb6b}
%C_3\big(|\nabla v_{r,R}|^{2-p} - |\nabla u|^{2-p}\big)^{-1}
%\leq I(p)\leq C_4\big(|\nabla v_{r,R}|^{2-p}+|\nablau|^{2-p}\big)^{-1},
%\end{equation}
%for $1<p<2$, with the same constants $C_3, C_4>0$. 
This may be written as
\begin{equation}
\label{nb6c}
C_9\big(|\nabla v_{r,R}|^{|p-2|} - |\nabla u|^{|p-2|}\big)
\leq I(p)\leq C_{10}\big(|\nabla v_{r,R}|^{|p-2|}+|\nabla
u|^{|p-2|}\big),
\end{equation}
where $C_9=\min\{C_7,1/C_8\}$ and $C_{10}=\max\{1/C_7,C_8\}$.

Thus by (\ref{nb5}), (\ref{nb6c}) we find,
\begin{multline}
\int_{O\cap O_2}
\varphi^2|\nabla v_{r,R}-\nabla u|^2
\big(|\nabla v_{r,R}|^{|p-2|}+|\nabla u|^{|p-2|}\big)dm\\
\leq C_{11}\int_{O\cap O_2}
|\varphi||v_{r,R}-u||\nabla \varphi||\nabla v_{r,R}-\nabla u|
\big(|\nabla v_{r,R}|^{|p-2|}+|\nabla u|^{|p-2|}\big)dm\\
\leq
C_{11}\bigg(\int_{O\cap O_2}
|\nabla\varphi|^2|v_{r,R}-u|^2
\big(|\nabla v_{r,R}|^{|p-2|}+|\nabla u|^{|p-2|}\big)dm
\bigg)^{1/2}\\
\cdot\bigg(
\int_{O\cap O_2}\varphi^2
|\nabla v_{r,R}-\nabla u|^2
\big(|\nabla v_{r,R}|^{|p-2|}+|\nabla u|^{|p-2|}\big)dm
\bigg)^{1/2}
\end{multline}
and
\begin{multline*}
\int_{O\cap O_2}\varphi^2
|\nabla v_{r,R}-\nabla u|^2
\big(|\nabla v_{r,R}|^{|p-2|}+|\nabla u|^{|p-2|}\big)dm\\
\leq
C_{11}^2\int_{O\cap O_2}|\nabla\varphi|^2
|v_{r,R}-u|^2
\big(|\nabla v_{r,R}|^{|p-2|}+|\nabla u|^{|p-2|}\big)dm.
\end{multline*}
Remembering (\ref{nb1}) we have
\begin{multline*}
\int_{O\cap O_1}|\nabla v_{r,R}-\nabla u|^2
\big(|\nabla v_{r,R}|^{|p-2|}+|\nabla u|^{|p-2|}\big)dm\\
\leq
C_{11}^2\int_{D_{r,R}\cap (O_2\setminus \overline{O_1})}
|\nabla\varphi|^2 |v_{r,R}-u|^2
\big(|\nabla v_{r,R}|^{|p-2|}+|\nabla u|^{|p-2|}\big)dm.
\end{multline*}
Because $\varphi$ is constant on $\Sigma_k(t)$ and $|\nabla
d_k|\equiv1$,
we have  by Theorem \ref{coarea}
\begin{multline*}
\int_{D_{r,R}\cap (O_2\setminus \overline{O_1})}
|\nabla\varphi|^2 |v_{r,R}-u|^2
\big(|\nabla v_{r,R}|^{|p-2|}+|\nabla u|^{|p-2|}\big)dm \\
\leq \int_{\{x : \varepsilon_1 < d_k(x) < \varepsilon_2\}}|\nabla\varphi|^2 |v_{r,R}-u|^2
\big(|\nabla v_{r,R}|^{|p-2|}+|\nabla u|^{|p-2|}\big)dm\\ =
\int_{\varepsilon_1}^{\varepsilon_2}
\eta'^{\, 2} H(t)dt,
\end{multline*}
where 
\begin{equation}
H(t)=\int_{\Sigma_k(t)}|v_{r,R} -u|^2
\big(|\nabla v_{r,R}|^{|p-2|}+|\nabla u|^{|p-2|}\big)
d{\Hm}^{n-1}.
\end{equation}
By H\"older's inequality
\begin{equation*}
1 \leq \int_{\varepsilon_1}^{\varepsilon_2}\eta'
H(t)^{1/2}H(t)^{-1/2} dt  \leq 
\bigg(\int_{\varepsilon_1}^{\varepsilon_2}\eta'^{\, 2} H(t)dt\bigg)^{1/2}
\bigg(\int_{\varepsilon_1}^{\varepsilon_2}H^{-1}(t)dt\bigg)^{1/2}.
\end{equation*}
It follows that
\begin{equation}
\label{ny1}
\bigg(\int_{\varepsilon_1}^{\varepsilon_2}H^{-1}(t)dt\bigg)^{-1}
\leq
\int_{\varepsilon_1}^{\varepsilon_2}\eta'^{\, 2}H(t)dt,
\end{equation}
for all $\varphi(x)=\eta(d_k(x))$ satisfying (\ref{nb1}).

We define a function $\hat \eta$ by the formula
$$
\hat \eta(s)=
\bigg(\int_{\varepsilon_1}^s H^{-1}(t) dt\bigg)
\bigg(\int_{\varepsilon_1}^{\varepsilon_2} H^{-1}(t) dt\bigg)^{-1}.
$$
Now $\hat \eta(\varepsilon_1)=0$ and $\hat\eta(\varepsilon_2)=1$.
Because
$$
\hat \eta'(s) =
\frac{1}{H(s)} \bigg(\int_{\varepsilon_1}^{\varepsilon_2}
\frac{dt}{H(t)}\bigg)^{-1},
$$
we have by (\ref{ny1})
\begin{multline*}
\bigg(\int_{\varepsilon_1}^{\varepsilon_2}H^{-1}(t)dt\bigg)^{-1}
\leq
\inf_{\varphi}\int_{\varepsilon_1}^{\varepsilon_2}\eta'^{\, 2}H(t)dt\\
\leq 
\int_{\varepsilon_1}^{\varepsilon_2}\hat\eta'^{\, 2}H(t)dt =
\bigg(\int_{\varepsilon_1}^{\varepsilon_2}H^{-1}(t)dt\bigg)^{-1}.
\end{multline*}
Because
\begin{multline*}
\int_{\varepsilon_1}^{\varepsilon_2}H^{-1}(t)dt\\=
\int_{\varepsilon_1}^{\varepsilon_2}dt \bigg(\int_{\Sigma_k(t)}|v_{r,R}
-u|^2
\big(|\nabla v_{r,R}|^{|p-2|}+|\nabla u|^{|p-2|}\big)
d{\Hm}^{n-1}\bigg)^{-1}\to \infty,
\end{multline*}
as $\varepsilon_2\to \infty$, the claim follows.
\qed

\section{Proofs of the corollaries}

\subsection*{Proof of Corollary \ref{cora}}

Let
\begin{equation}
H(t)=
\int_{\Sigma_k(t)}|v_{r,R} -u|^2
\big(|\nabla v_{r,R}|^{|p-2|}+|\nabla u|^{|p-2|}\big)
d{\Hm}^{n-1}.
\end{equation}
By H\"older's inequality
\begin{multline*}
(S-r)^2 = \Big(\int_r^Sdt\Big)^2
= \Big(\int_r^S\frac{H^{-1/2}(t)}{H^{-1/2}(t)}dt\Big)^2
 \leq
\Big(\int_r^SH^{-1}(t)dt\Big)
\Big(\int_r^S H(t) dt\Big).
\end{multline*}
Hence
\begin{equation}
\label{hest}
(S-r)^2\Big(\int_r^SH^{-1}(t)dt\Big)^{-1}
\leq
\Big(\int_r^S H(t) dt\Big).
\end{equation}
Now by (\ref{hest}) and Theorem \ref{coarea}
\begin{multline*}
\Bigg[\int_r^S dt
\bigg(\int_{\Sigma_k(t)}|v_{r,R} -u|^2
\big(|\nabla v_{r,R}|^{|p-2|}+|\nabla u|^{|p-2|}\big)
d{\Hm}^{n-1}\bigg)^{-1}\Bigg]^{-1}\\
\leq
\frac{1}{(S-r)^2}
\int_{D_{r,S}}|v_{r,R} -u|^2
\big(|\nabla v_{r,R}|^{|p-2|}+|\nabla u|^{|p-2|}\big)
dm\to 0,
\end{multline*}
as $S\to \infty$, proving the claim.
\qed

\subsection*{Proof of Corollary \ref{corb}}
Since
$$
\int_{D_{r,\infty}}|v_{r,R} -u|^2
\big(|\nabla v_{r,R}|^{|p-2|}+|\nabla u|^{|p-2|}\big)
dm \leq M < \infty,
$$
we have for $S>r$, 
$$
\frac{1}{S^2} \int_{D_{r,S}}|v_{r,R} -u|^2
\big(|\nabla v_{r,R}|^{|p-2|}+|\nabla u|^{|p-2|}\big)
dm \leq \frac{M}{S^2}\to 0,
$$
as $S\to\infty$.
\qed

%\subsection*{Acknowledgements}

%\newpage

%\bigskip

\end{document}